\documentclass[12 pt]{article}
\usepackage[dvips]{graphicx}
\usepackage{amssymb}

\newcommand \no{\noindent}
\newcommand \implies{\Rightarrow}

\newcommand{\be}{\begin{equation}}
\newcommand{\ee}{\end{equation}}
\newcommand{\bea}{\begin{eqnarray}}
\newcommand{\eea}{\end{eqnarray}}
\def\reff#1{(\ref{#1})}

\newtheorem{theorem}{Theorem}

\begin{document}
\bibstyle{ams}

\title{Compact packings of the plane with two sizes of discs}

\author{Tom Kennedy
\\Department of Mathematics
\\University of Arizona
\\Tucson, AZ 85721
\\ email: tgk@math.arizona.edu
\bigskip
}

\maketitle

\begin{abstract} We consider packings of the plane using discs of 
radius $1$ and $r$. A packing is compact if every disc $D$ is 
tangent to a sequence of discs $D_1, D_2, \cdots, D_n$ such that 
$D_i$ is tangent to $D_{i+1}$. We prove that there are only nine values of 
$r$ with $r<1$ for which such packings are possible. For each of the nine
values we describe the possible compact packings.
\end{abstract}


\vspace{\fill}
\hrule width2truein
\smallskip
{\baselineskip=12pt
\noindent
\copyright\ 2004 by the author. Reproduction of
this article is permitted for non-commercial purposes.
\par}

\newpage

\section{Introduction} \setcounter{equation}{0}

\bigskip

A packing of the plane with discs is said to be compact if 
every disc $D$ is tangent to a sequence of discs 
$D_1, D_2, \cdots, D_n$ such that $D_i$ is tangent to $D_{i+1}$
for $i=1,2,\cdots,n$ with $D_{n+1}=D_1$ \cite{fta}. 
If we pack the plane using discs of the same
radius, then the only possible compact packing is the triangular 
lattice in which each disc is surrounded by six discs tangent to it.  
In this paper we ask what compact packings are possible if we pack the 
plane using discs of radius $1$ and $r$ with $r<1$. We do not impose any 
condition on the relative number of discs of the two radii except
to require that discs of both sizes be present. 
We will prove that there are only nine values of $r$ for which a 
compact packing is possible. For these values we describe
the possible compact packings. 

\begin{figure}[tbh]
\includegraphics{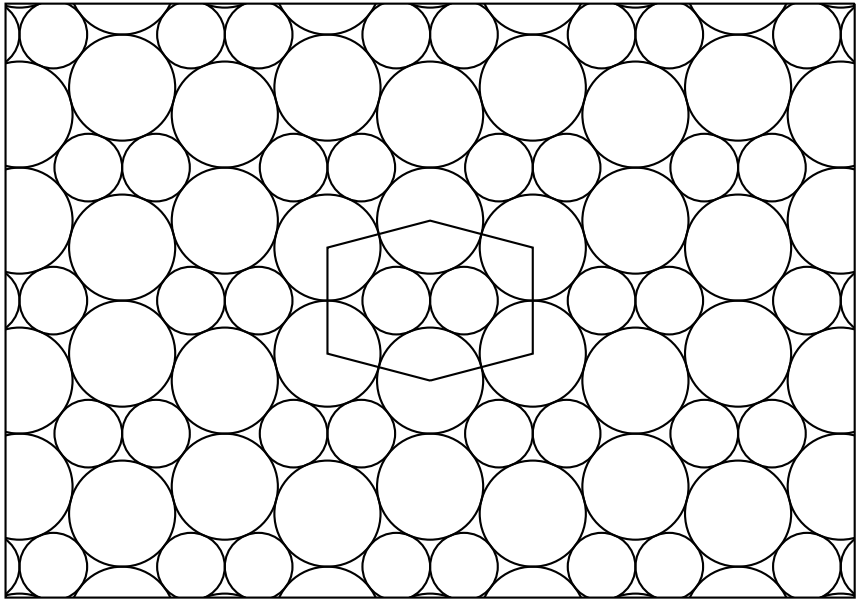}
\caption{A compact packing with $r=c_1=0.6375559772 \cdots$. $r$ 
is a root of  $r^4 - 10r^2 - 8r + 9$. }
\label{fig637}
\end{figure}

\begin{figure}[tbh]
\includegraphics{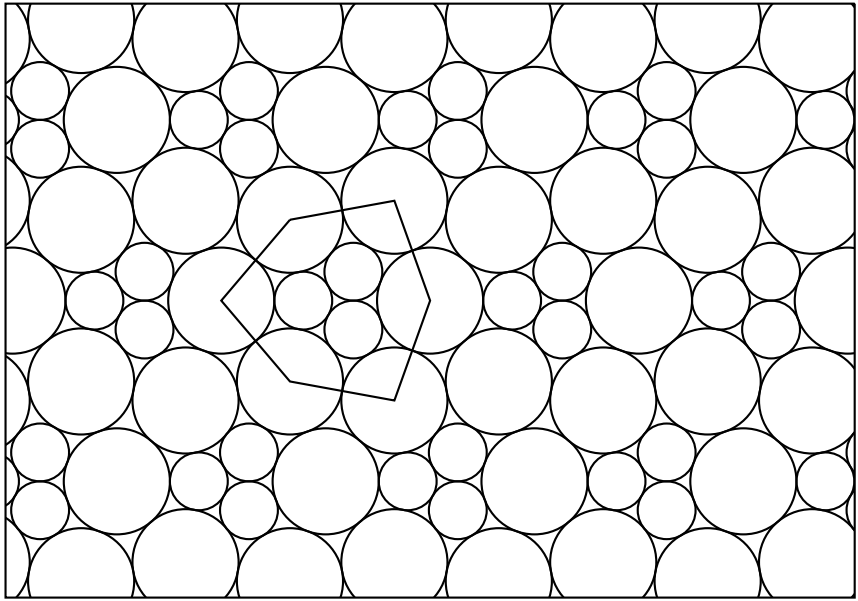}
\caption{A compact packing with $r=c_2=0.5451510421 \cdots$. 
$r$ is a root of 
$(7+4\sqrt{3}) r^4 + (20+12\sqrt{3}) r^3 
+ (6+4\sqrt{3}) r^2+(-20-4\sqrt{3})r+3$.
}
\label{fig545}
\end{figure}

\begin{figure}[tbh]
\includegraphics{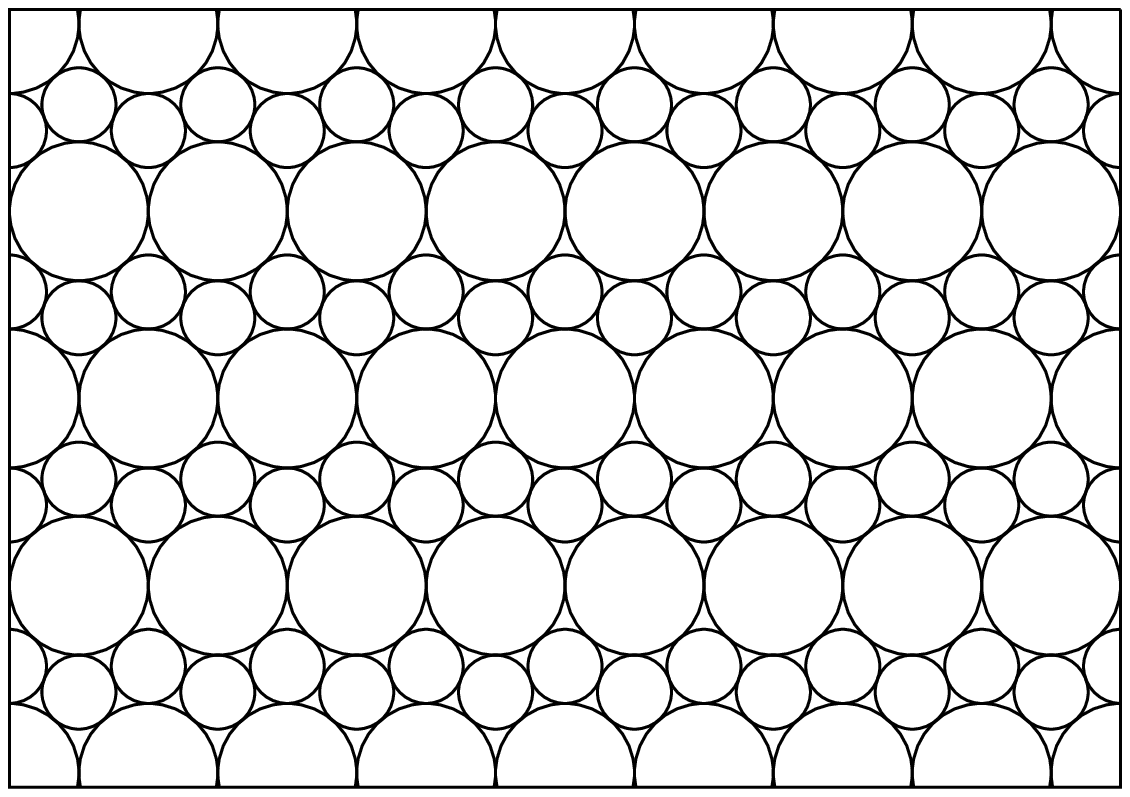}
\caption{A compact packing with $r=c_3=0.5332964167 \cdots$. 
$r$ is a root of $8r^3 + 3r^2 - 2r - 1$. 
}
\label{fig533}
\end{figure}

\begin{figure}[tbh]
\includegraphics{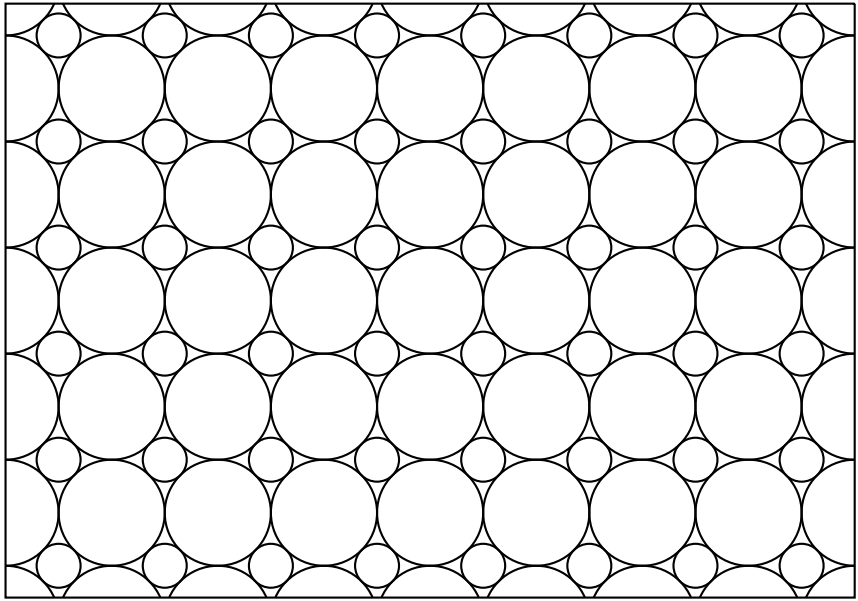}
\caption{A compact packing with $r=c_4=0.4142135624 \cdots =\sqrt{2}-1$.
}
\label{fig414}
\end{figure}

\begin{figure}[tbh]
\includegraphics{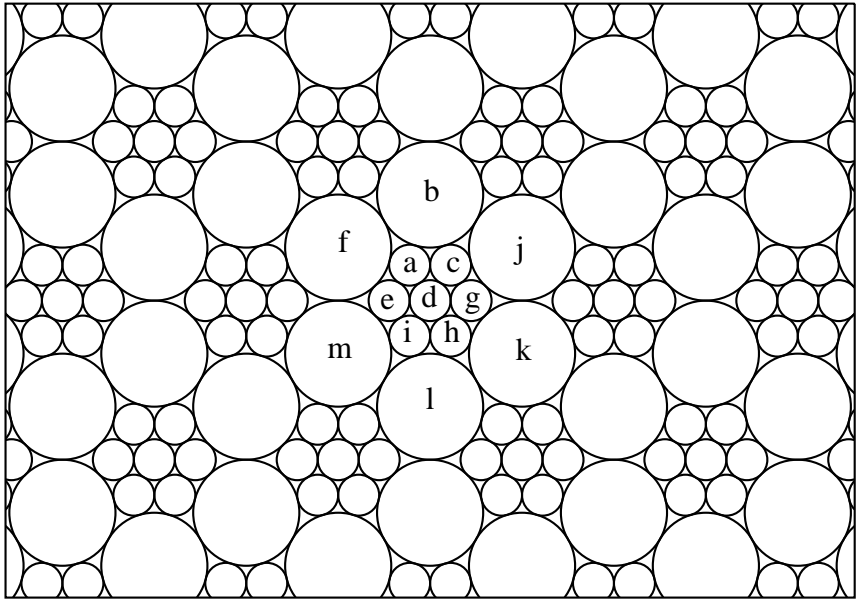}
\caption{A compact packing with $r=c_5=[2 \sqrt{3}+1 -2 \sqrt{1+\sqrt{3}}]/3
=0.3861061049 \cdots$. 
}
\label{fig386}
\end{figure}

\begin{figure}[tbh]
\includegraphics{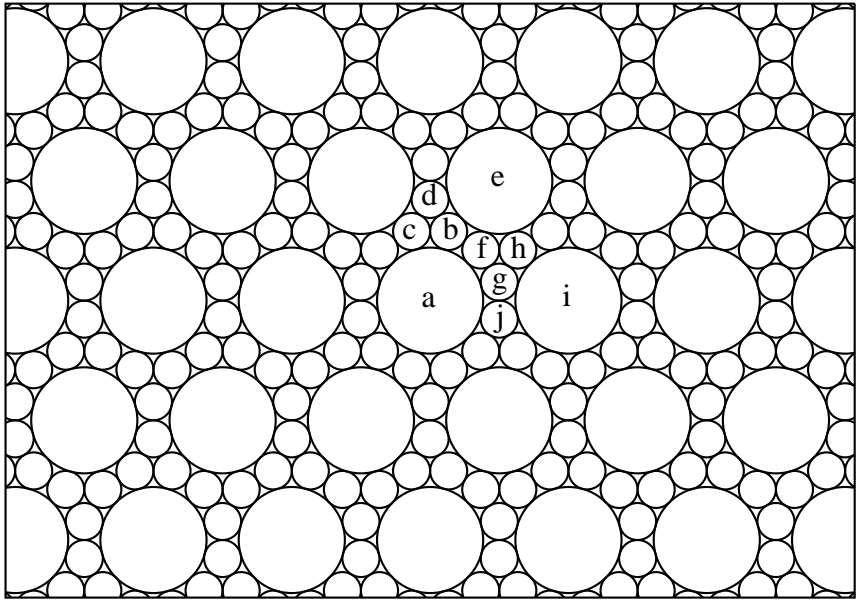}
\caption{The compact packing with $r=c_6=\sin(\pi/12)/(1-\sin(\pi/12))
=0.3491981862 \cdots$. 
}
\label{fig349}
\end{figure}

\begin{figure}[tbh]
\includegraphics{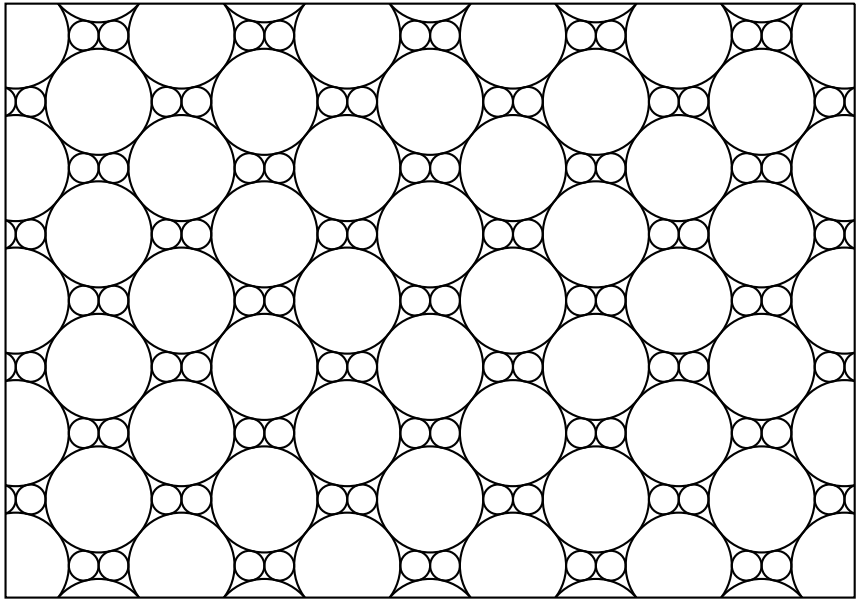}
\caption{A compact packing with $r=c_7=(\sqrt{17}-3)/4
=0.2807764064 \cdots$. 
}
\label{fig280}
\end{figure}

\begin{figure}[tbh]
\includegraphics{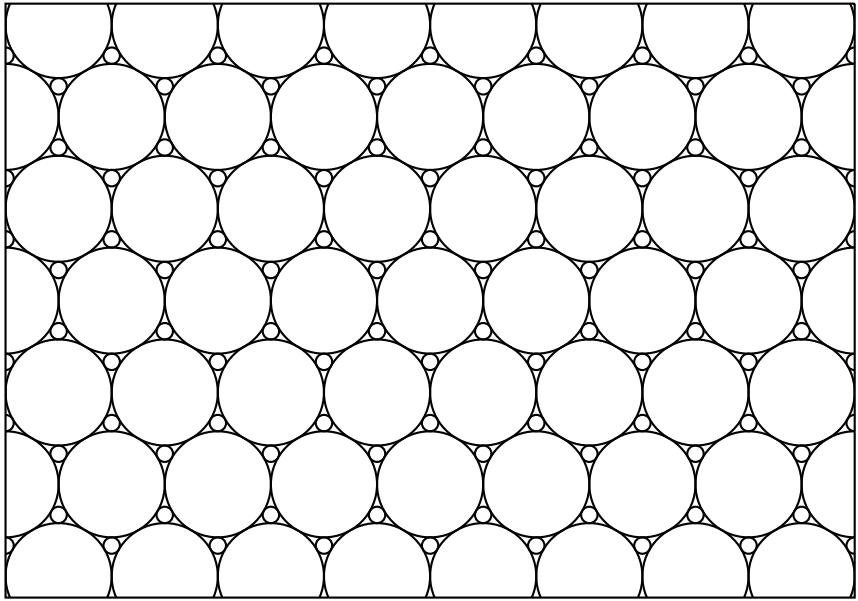}
\caption{A compact packing with $r=c_8= 2 \sqrt{3}/3-1 = 0.1547005384 \cdots$.
}
\label{fig154}
\end{figure}

\begin{figure}[tbh]
\includegraphics{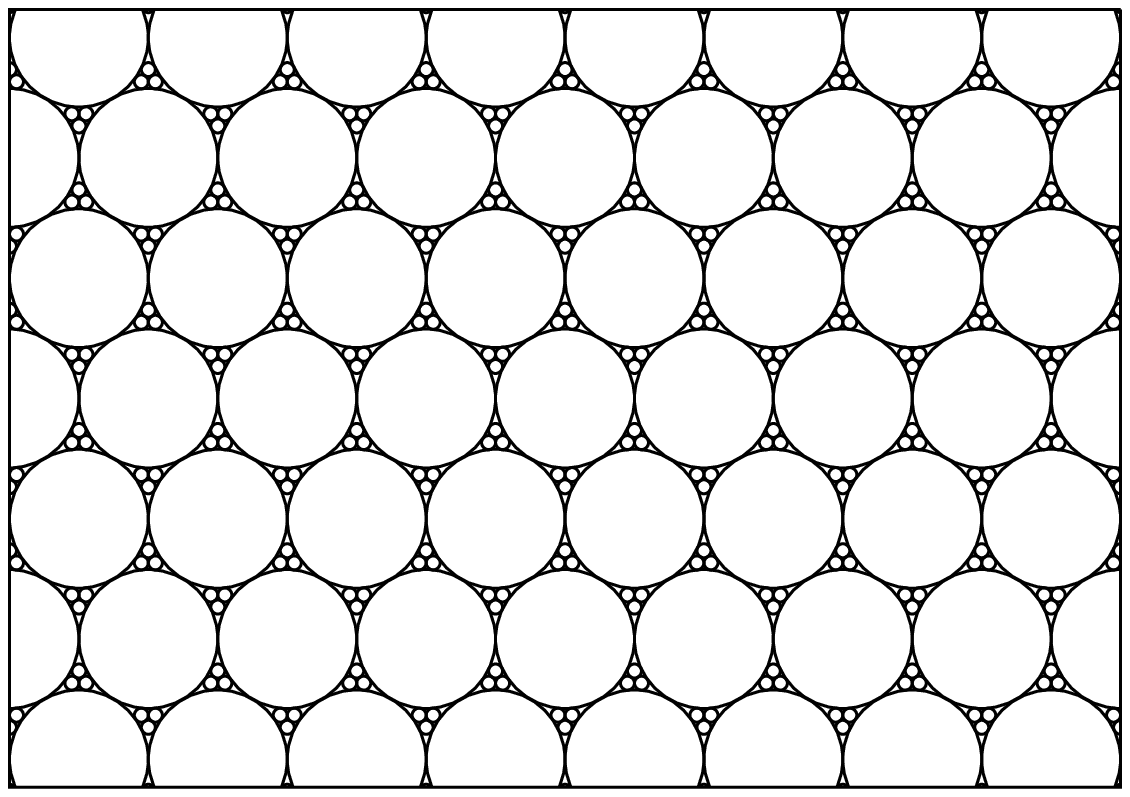}
\caption{A compact packing with $r=c_9=5-2 \sqrt{6}=0.1010205144 \cdots$.  
}
\label{fig101}
\end{figure}

The nine values of $r$ which allow compact packings are given in 
table \ref{table_r}. Examples of compact packings for each of the 
nine values are given in figures \ref{fig637} to \ref{fig101}.
Compact packings for seven of the nine possible values of $r$ 
(all but $c_2$ and $c_5$) appear 
in L. Fejes T\'oth's classic book \cite{ftb} and in \cite{hm,mol}.
A packing with $r=c_5$ appears in \cite{lh}. 

\begin{table}[thp]
\begin{center}
\begin{tabular}{ | r | c | c | c | c | }         \hline
  & decimal & exact & i j k & sequence \\ \hline \hline
 $c_1$ & 0.6375559772 & $r^4 - 10r^2 - 8r + 9=0$  & 3 2 0 & 1111r \\ \hline
 $c_2$ & 0.5451510421 & $P(r)=0$  & 2 2 1 & 111rr \\ \hline
 $c_3$ & 0.5332964167 & $8r^3 + 3r^2 - 2r - 1=0$  & 1 4 0 &  1r1r1 \\ \hline
 $c_4$ & 0.4142135624 & $\sqrt{2}-1$  & 4 0 0 & 1111 \\ \hline
 $c_5$ & 0.3861061049 & $[2 \sqrt{3}+1 -2 \sqrt{1+\sqrt{3}}]/3$ & 1 2 2 & 
   1rrr1 \\ \hline
 $c_6$ & 0.3491981862 & $\sin(\pi/12)/(1-\sin(\pi/12))$  & 0 4 1 & 1rr1r  
   \\ \hline
 $c_7$ & 0.2807764064 & $(\sqrt{17}-3)/4$  & 2 2 0 & 111r \\ \hline
 $c_8$ & 0.1547005384 & $2 \sqrt{3}/3-1$  & 3 0 0 & 111 \\ \hline
 $c_9$ & 0.1010205144 & $5-2 \sqrt{6} $  & 1 2 1 & 11rr \\ \hline
\end{tabular}
\caption
 {\protect The nine values of $r$ for which a compact packing using discs
   of radius $1$ and $r$ is possible. The decimal expressions are 
   of course approximations. For six values of $r$ we give 
   an explicit exact expression. For the other three we give polynomials
   which have $c_i$ as a root. In particular, $c_2$ is a root of 
   $P(r)=(7+4\sqrt{3}) r^4 + (20+12\sqrt{3}) r^3 + (6+4\sqrt{3}) 
      r^2+(-20-4\sqrt{3})r+3$. 
   The last two columns give information on the sequence of discs that 
   can appear around a small disc. These two columns are explained in 
   later sections.
 }
\label{table_r}
\end{center}
\end{table} 

An interesting question is to find (for each $r$) the densest packing 
of the plane using discs of radius $1$ and $r$. 
For six of the above values of $r$ ($c_1, c_3, c_4, c_6, c_7, c_8$ ),
Heppes has proved that the densest packing is a compact packing 
\cite{hepa,hepb}. We expect that for the other three values of $r$ 
the densest packing is also a compact packing. Rigorous bounds on the 
densest packing may be found in \cite{bb,flor}, and a non-rigorous 
study is in \cite{lh}.

\begin{figure}[tbh]
\includegraphics{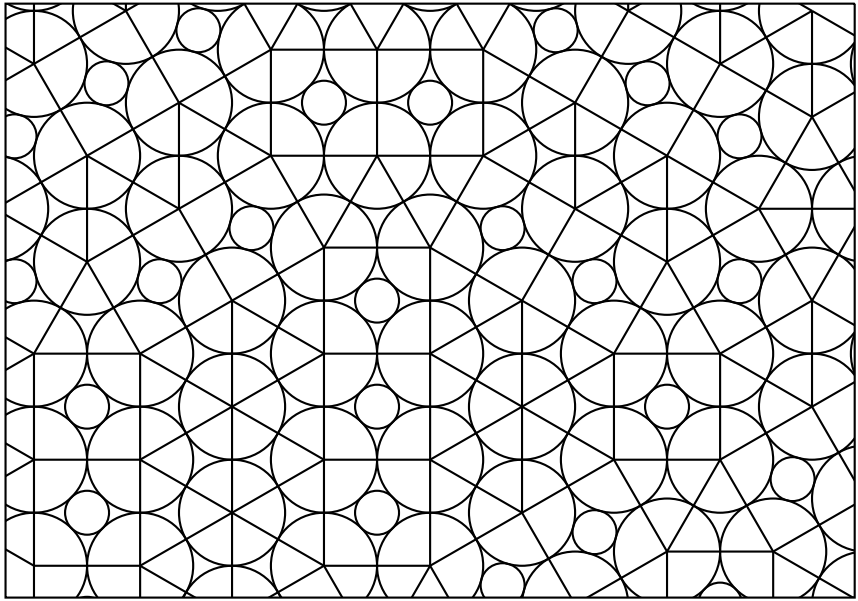}
\caption{A compact packing with $r=c_4$ and the corresponding
tiling by squares and triangles.
}
\label{fig414a}
\end{figure}

\begin{figure}[tbh]
\includegraphics{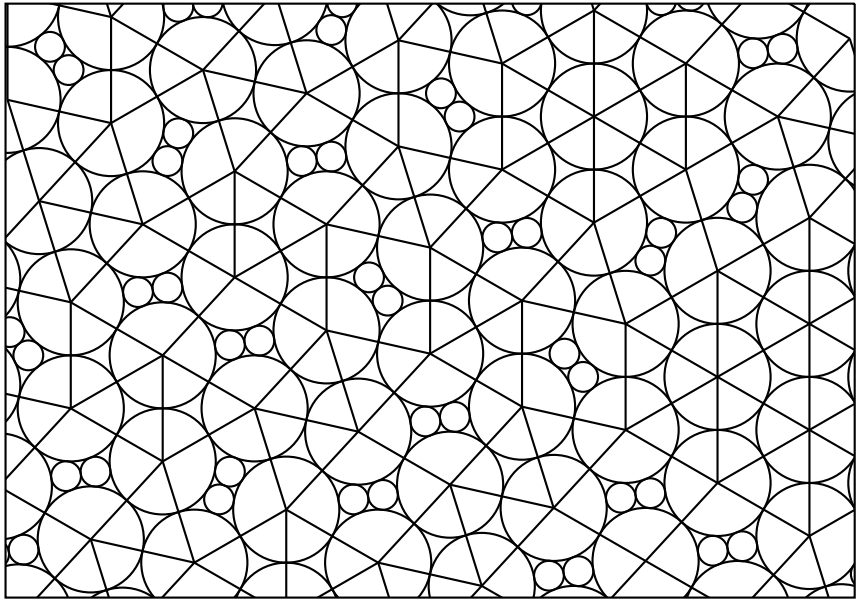}
\caption{A compact packing with $r=c_7$ 
and the corresponding tiling by rhombi and triangles.}
\label{fig280a}
\end{figure}

\begin{figure}[tbh]
\includegraphics{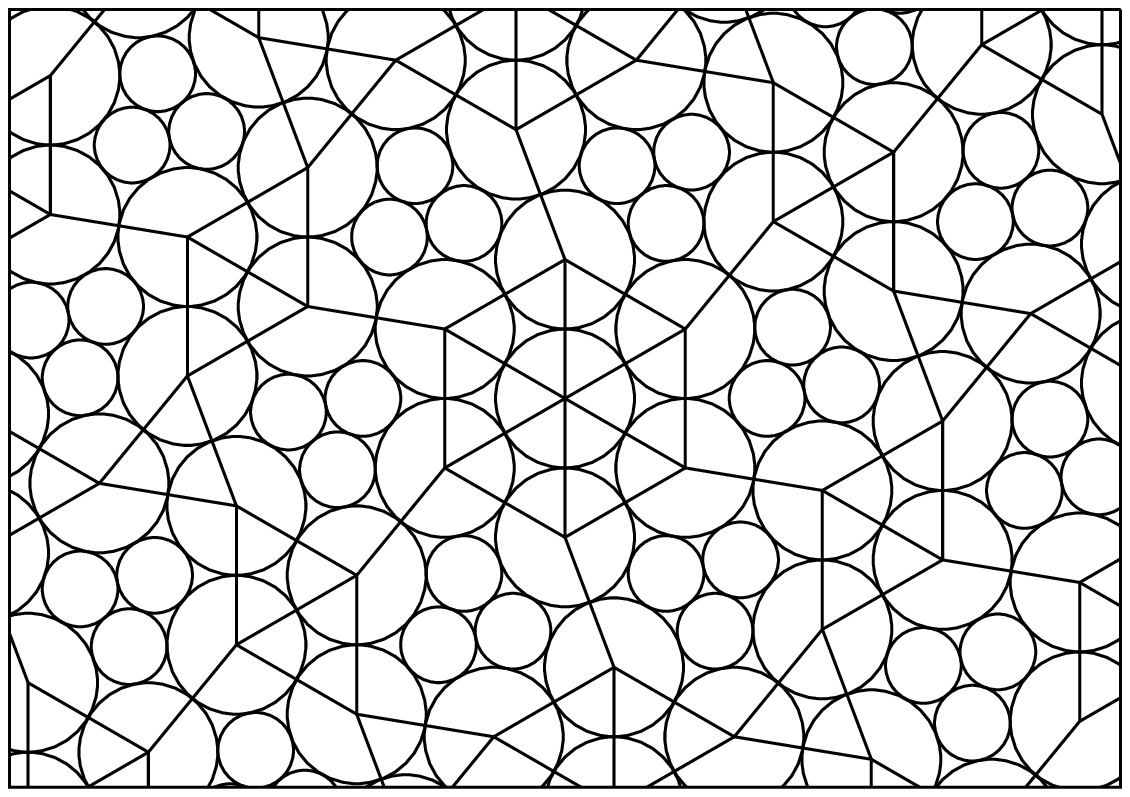}
\caption{Another compact packing with $r=c_2$.
} 
\label{fig545t}
\end{figure}

\begin{figure}[tbh]
\includegraphics{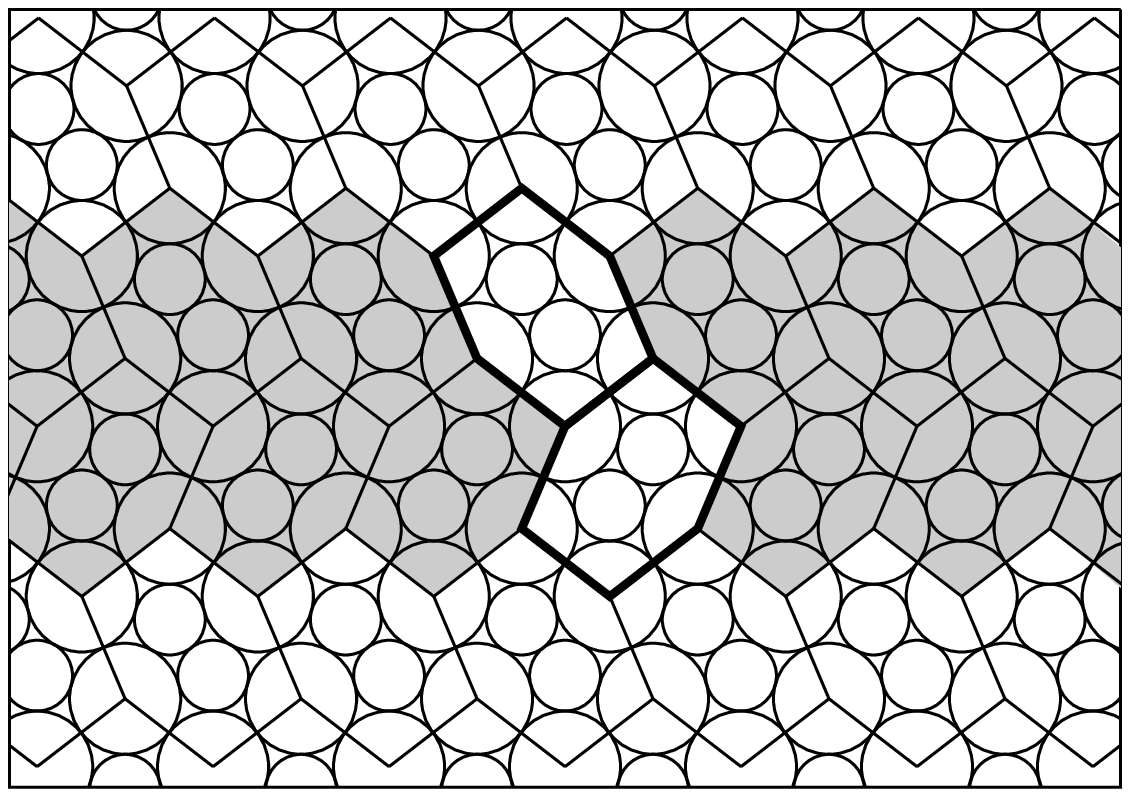}
\caption{Another compact packing with $r=c_1$.
}
\label{fig637a}
\end{figure}

\begin{figure}[tbh]
\includegraphics{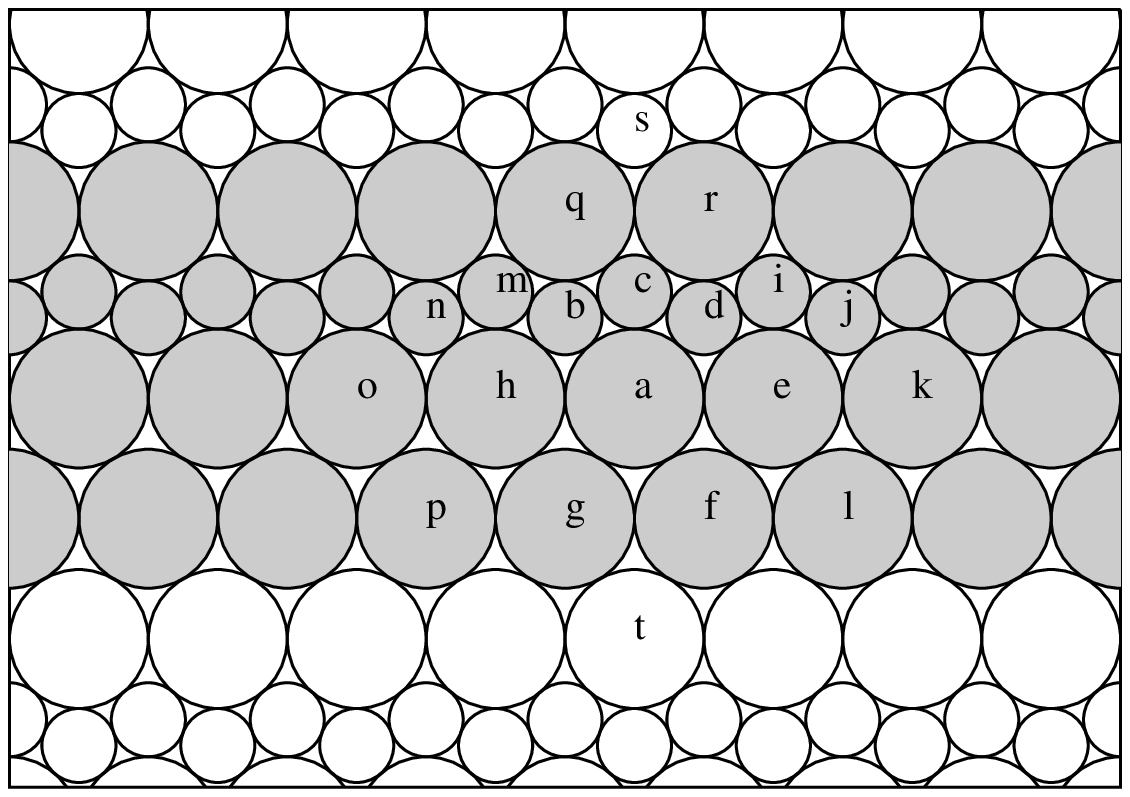}
\caption{Another compact packing with $r=c_3$. 
}
\label{fig533a}
\end{figure}

\begin{figure}[tbh]
\includegraphics{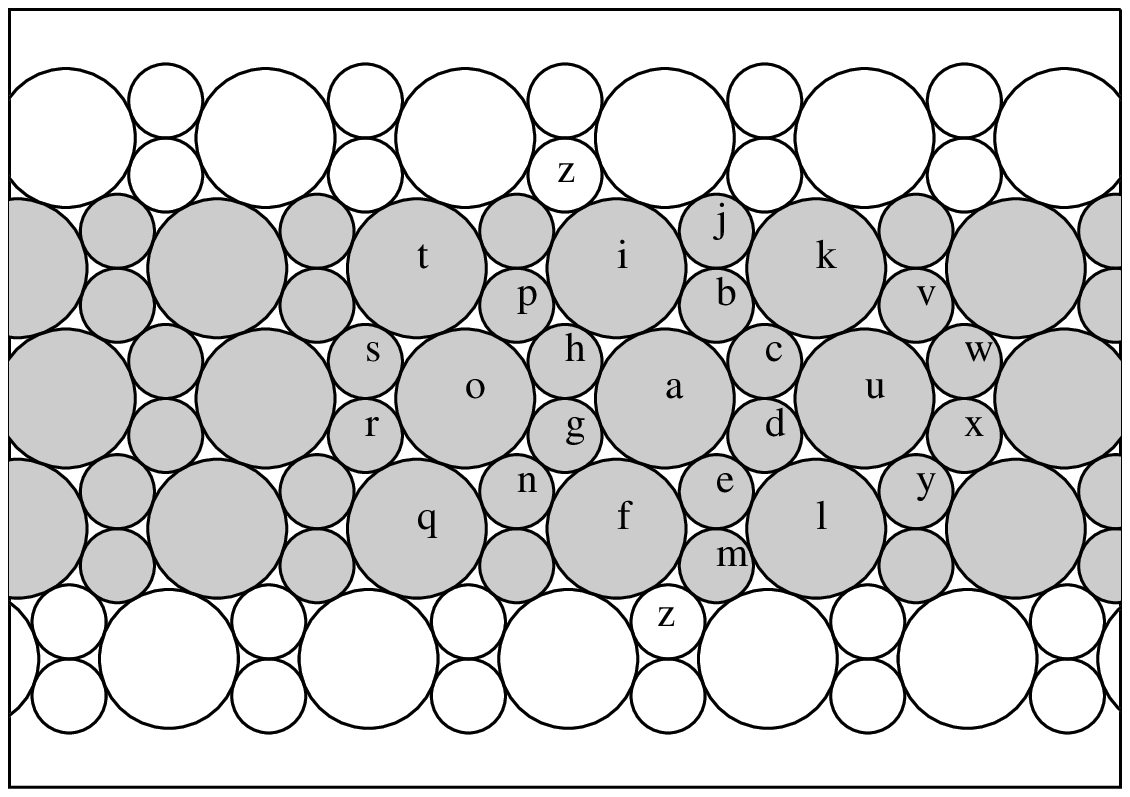}
\caption{Yet another compact packing with $r=c_3$. 
}
\label{fig533b}
\end{figure}

\section{Values of $r$ that allow compact packings}
\label{values}

\begin{theorem} The nine values of $r$ given in table \ref{table_r}
are the only values in $(0,1)$ for which there exists a compact packing 
of the plane using discs of radius $r$ and $1$. 
\end{theorem}

This section is devoted to the proof of this theorem. 
The idea is simple. Consider the center of some disc and the 
angles that are formed by drawing line segments from the center of the 
disc to the centers of the discs that are tangent to the original disc. 
There are only a few possible values for these angles and they are functions
of $r$. The sum of the angles around the center must be $2 \pi$, and 
this constrains $r$. 

To begin the proof, consider three discs that are tangent to each other
and the angles in the triangle whose vertices are their centers. 
If all three discs have radius $1$ or all three have radius
$r$, then the angles are all $\pi/3$. 
If two of the discs have radius $1$ and one has radius $r$, 
then the triangle has one angle of $\alpha$ and two of $\alpha^\prime$ where 
\be
\cos(\alpha^\prime) = {1 \over 1+r}, \quad \alpha=\pi- 2 \alpha^\prime
\label{alpha}
\ee
If two of the discs have radius $r$ and one has radius $1$, 
then the triangle has one angle of $\beta$ and two of $\beta^\prime$ where 
\be
\cos(\beta^\prime) = { r \over 1+r}, 
\quad \beta =\pi - 2 \beta^\prime
\label{beta}
\ee
Thus in a compact packing the angles around a disc of radius $r$ 
can only be $\pi/3$, $\alpha$ and $\beta^\prime$. 
And the angles around a disc of radius $1$ can only be 
$\pi/3$, $\alpha^\prime$ and $\beta$. 
Since the angles around a disc must sum to $2 \pi$, there must exist 
non-negative integers $i,j,k$ such that 
\be
i \alpha + j \beta^\prime + k \pi/3=2 \pi 
\label{coma}
\ee
and non-negative integers $l,m,n$ such that 
\be
l \alpha^\prime + m \beta + n \pi/3=2 \pi 
\label{comb}
\ee
For every value of $r$ there is a trivial solution of both of these equations,
namely $i=j=l=m=0, k=n=6$. 
In a compact packing that contains discs of both radii, there 
must be at least one pair of discs of different radii that are tangent,
and so there must be at least one solution of \reff{coma} other than the 
trivial one and at least one solution of \reff{comb} other than the trivial
one. We start by determining when \reff{coma} has solutions.

Define
\be
F_{ijk}(r)=i \alpha + j \beta^\prime + k \pi/3
\label{comf}
\ee
Equations \reff{alpha} and \reff{beta} imply that $\alpha$ and 
$\beta^\prime$ are decreasing functions of $r$, and so $F_{ijk}(r)$ 
is also decreasing for each choice of $ijk$. This implies that for 
each choice of $ijk$ there is at most one value of $r$ for which 
$F_{ijk}(r)=2\pi$. 

As $r$ goes to zero, 
$\alpha$ converges to $\pi$ and $\beta^\prime$ converges to $\pi/2$. 
So 
\be
\lim_{r \rightarrow 0} F_{ijk}(r) = i \pi + j \pi/2 + k \pi/3 
\ee
And when $r=1$ all the angles are $\pi/3$, so 
\be
F_{ijk}(1) =(i + j + k) \pi/3
\ee
So there is an $r \in (0,1)$  with $F_{ijk}(r)=2 \pi$ if and only if 
we have  $6i + 3j + 2k>12$ and $i + j + k<6$.
Thus there are only a finite number of $i,j,k$ for which eq. \reff{coma}
might have a solution. 
There are further constraints on $i,j,k$. Recall that the sides adjacent
to the angle $\alpha$ both have length $1+r$, the sides adjacent to 
$\beta$ have lengths $1+r$ and $2r$ and the sides adjacent to $\pi/3$ 
both have length $2r$. This implies that $j$ must be even, and 
if $j$ equals zero, then one of $i$ and $k$ must be zero. 
It is now trivial to check that the only possibilities for $ijk$ 
are the nine listed in the next to last column of the table and one 
additional case, $i=5, j=0, k=0$. For this last case, 
$F_{ijk}(r)=2 \pi$ at $r=(1-\sin(\pi/5))/\sin(\pi/5)=0.7013025284\cdots$.

The case of $r=(1-\sin(\pi/5))/\sin(\pi/5)$ can be ruled out by 
considering eq. \reff{comb}. We compute $\alpha^\prime$ and $\beta$ 
and check if \reff{comb} is satisfied for any choice of $l,m,n$. 
We find that it is not. Thus the only $r$ which might allow compact
packings are those shown in the table. The explicit packings 
shown in figure \ref{fig637} to \ref{fig101} show that compact packings 
with these nine values of $r$ are indeed possible. 
It is straightforward to use the figures to compute the nine values of $r$ 
given in the table.
This completes the proof of theorem 1. 

\section{The compact packings}
\label{packings}

In this section we will describe the possible compact packings 
for $r=c_1$ to $c_9$. 
It will be helpful to put eqs. \reff{coma} and \reff{comb} in a slightly
different form. Consider the triangle formed by the 
centers of discs with radii $r_1,r_2,r_3$ which are tangent to each other.
Let $\theta(r_1,r_2,r_3)$ denote the angle at the center of the 
disc of radius $r_1$. Depending on the values of $r_1,r_2,r_3$, 
$\theta(r_1,r_2,r_3)$ is either 
$\pi/3,\alpha,\alpha^\prime,\beta$ or $\beta^\prime$.
Now let $D$ be a disc with radius $r_0$ and let 
$D_1,D_2,\cdots D_n$ be the discs that are tangent to $D$ and such that 
$D_i$ is tangent to $D_{i+1}$. Let $r_1,r_2,\cdots,r_n$ be their radii.
Then we must have 
\be
\sum_{i=1}^n \theta(r_0,r_i,r_{i+1}) = 2 \pi 
\label{seqeq}
\ee
where $r_{n+1}=r_1$. Of course each $r_i$ is either $1$ or $r$. 
So for $r_0$ equal to $1$ or $r$, we want to find all finite sequences 
$r_1,r_2,\cdots,r_n$ consisting of $1$'s and $r$'s such that the 
above equation is satisfied. (As we noted before, we seek solutions beyond the
trivial ones with $n=6$ and $r_1=\cdots=r_6=r_0$ which are solutions
for all $r$.)

We will write the cases for which \reff{seqeq} is satisfied in the form 
\be 
r_0 : r_1 r_2 r_3 \cdots r_n
\label{rseq}
\ee
meaning that a disc of radius $r_0$ can be surrounded by a sequence
of discs of radius $r_1,r_2,\cdots,r_n$. Of course, if  
$r_0 : r_1 r_2 r_3 \cdots r_n$ is possible then so is 
$r_0 : r_n r_{n-1} r_{n-2} \cdots r_1$ as well as any 
cyclic permutation of $r_1 r_2 r_3 \cdots r_n$. We will only list one
possibility from each equivalence class. 
We will refer to \reff{rseq} as a ``sequence for small discs'' when 
$r_0=r$ and as a ``sequence for large discs'' when $r_0=1$.
When there is only one possible sequence up to the symmetries just 
described, we will say the sequence is unique. 

We have already determined the solutions of eq. \reff{coma} and thus the 
solutions of \reff{seqeq} with $r_0=r$. 
For each of the nine values of $r$ there are two possible sequences,
$r:rrrrrr$ and the sequence given in the last column of the table.
We will argue that $r:rrrrrr$ does not occur 
except when $r=c_5$. Suppose the configuration contains at least one
small disc with sequence $rrrrrr$. Consider the small discs that are 
tangent to at least one small disc with sequence $rrrrrr$. They cannot
all have sequence $rrrrrr$. So there is at least one small disc $D$ 
whose sequence is not $rrrrrr$ but which is tangent to a small disc with 
sequence $rrrrrr$. This implies that the sequence for $D$ must contain 
$rrr$ but not be equal to $rrrrrr$. So it must contain $1rrr$ or $rrr1$. 
From the table we see that this is possible only for $r=c_5$. 

Finding the possible solutions of \reff{seqeq} 
with $r_0=1$ is straightforward.
For each of the nine values of $r$ we compute 
$\alpha^\prime$ and $\beta$ and check for what choices of 
$r_1,r_2,\cdots,r_n$ eq. \reff{seqeq} is satisfied with $r_0=1$. We will 
list the possibilities as we consider each of the nine values of $r$. 

For some values of $r$ we can give very explicit descriptions of
all possible compact packings, but for others we can only show that 
there is a correspondence between packings and tilings using certain
polygons. 
The order in which we will consider the nine values of $r$ is
for the sake of exposition.
We start with the only value of $r$ for which the packing is unique. 
Of course, we can translate, rotate and reflect any compact packing
to get another compact packing. When we say unique we implicitly mean
unique up to these symmetries. Throughout the arguments that follow
we will freely use these symmetries to assume, for example, that a certain
feature of the configuration is oriented in a particular direction. 

For the next three values of $r$, the compact configurations are obtained 
by starting with a triangular lattice of large discs and either 
replacing some of the discs by groups of small discs ($r=c_5$) or 
filling in some of the holes with small discs ($r=c_8,c_9$).
For the next three values, $c_4,c_7,c_2$, we can only show that there 
is a one to one correspondence between compact packings and tilings 
of the plane with polygons whose sides all have length $2$. 
For the last two values, $c_1,c_3$, we show that the packing is made of 
layers, but for each layer there is a two-fold choice.

\medskip

\no {\bf Possible packings for $r=c_6$:} 

There is a unique sequence for small discs and apriori two possible sequences 
for large discs. 
\bea
&r : 1rr1r \nonumber\\
&1 : rrrrrrrrrrrr \quad or \quad 111111
\eea
We will argue that the $1:111111$ sequence does not appear. The argument 
is similar to the argument used to rule out $r:rrrrrr$. 
If $1:111111$ appears, then there must be a large disc $D$ whose sequence 
is not $111111$ but which is tangent to a large disc with sequence 
$111111$. So the sequence of $D$ must contain $111$ but not equal $111111$. 
This contradicts the above possible sequences around a large disc. 
Thus the only possible sequence around a large disc is $rrrrrrrrrrr$. 

We will show that these unique sequences imply the configuration is unique. 
We start with a large disc which we label $a$ in figure \ref{fig349}. 
By the sequence for large discs it is surrounded by 12 small discs. 
Consider the small disc labeled $b$ in the figure. 
Up to symmetry, there is only one sequence
that can surround it. Reflecting the configuration if necessary, we
can assume that the discs around disc $b$ are as shown in the figure, i.e.,
discs $a,c,d,e,f$. Applying the sequence for small discs to disc $f$, we
conclude that discs $g$ and $h$ must be present. Applying the sequence 
for small discs to disc $g$, we conclude that discs $i$ and $j$ are present. 
Note that we now have three large discs, $a$, $e$ and $i$, forming 
an equilateral triangle, and by the sequence for large discs each of them is 
surrounded by $12$ small discs. We can now repeat this argument
to fill out the configuration. 

\medskip

\no {\bf Possible packings for $r=c_5$:} 

This is the case for which there are two possible sequences for small discs.
For the possible sequences for a large disc we give 
a number after the sequence to indicate the number of discs around 
the large disc. (This is just to aid in reading the following.)
\bea
& r : 1rrr1 \quad or \quad rrrrrr \nonumber\\
& 1 : rr11111(7), \quad rrr1r111(8), \quad rrr11r11(8), 
\quad rr1rr111(8), \quad \nonumber\\
& rr11rr11(8), \quad rrrr1r1r1(9), \quad rrr1rr1r1(9), 
\quad rr1rr1rr1(9), \quad or \quad 111111(6)
\eea

In the packing shown in figure \ref{fig386} we can replace any hexagonal
block of 7 small discs by a single large disc. If we replace all such
blocks of 7 small discs, then
the result will just be the triangular packing of large discs. 
If we start with the triangular packing by large discs, 
take a subset of the discs such that no two of them are tangent, 
and replace the discs in this subset by hexagonal blocks of 7 small discs,
then the resulting packing will be compact. 
We claim that all compact packings can be obtained in this fashion.

We begin by showing that every connected component of small discs 
must consist of seven small discs arranged in a hexagon like those in 
figure \ref{fig386}.
A connected component of small discs must contain a small disc that is 
tangent to a large disc. We label this disc $a$ in figure \ref{fig386}.
Since it is tangent to a large disc, there is up to symmetry
only one possible sequence around it. We can assume the sequence is 
that shown in the figure, i.e., discs $b,c,d,e,f$.  
The possible sequences around $c,d$ and $e$ then imply that discs 
$g,h,i,j,m$ must be present. The possible sequences around $h$ then 
imply discs $k,l$ must be present. Thus small discs can only appear as 
part of a cluster of seven small discs arranged in a hexagon and 
surrounded by six large discs. 

Now we replace every hexagonal block of 7 small discs 
in our compact configuration by a single large disc. 
The result is a compact configuration with only large discs.
The only such configuration is the triangular packing. Thus our 
original packing is obtained by replacing some subset of non-tangent
large discs in the triangular packing by large discs by hexagonal 
blocks of 7 small discs. 

\medskip

\no {\bf Possible packings for $r=c_8$:} 

There is a unique sequence for small discs, 
\be
  r : 111,
\ee
and 13 sequences for large discs which we do not list. 
They are not needed in the argument.

The unique sequence for small discs says that each small disc is 
surrounded by three large discs that are tangent to each other. 
So if we remove the small discs, the configuration will still be compact. 
Thus we can remove all the small discs from our configuration and obtain 
a compact configuration with only large discs. There is only one such 
configuration, the triangular one. Thus we have shown that every 
compact configuration is given by putting small discs into some 
subset of the holes in the triangular lattice configuration of large discs. 
(Equivalently, every compact configuration can be obtained by removing 
some subset of the small discs from the configuration in figure
\ref{fig154}.)

\medskip

\no {\bf Possible packings for $r=c_9$:} 

The sequence for small discs is unique: 
\be
r : 11rr 
\ee
We do not list the $195$ possible sequences for large discs. 
Note that in the configuration in figure \ref{fig101} we have a triangular 
lattice of large discs with each hole filled with a cluster of 
three small discs which are tangent to each other. We can remove any 
of the clusters of three tangent small discs and still have a 
compact configuration. We will show that any compact configuration 
can be obtained in this way. 

Start with a small disc. The sequence for small discs implies it is tangent 
to two other small discs which are tangent to each other. 
The sequence for small discs applied to these three small discs then implies
that they sit inside the hole formed by three large discs that are 
tangent to each other. We can remove the three small discs and still 
have a compact configuration. Thus we can remove all the small discs 
and have a compact configuration with only large discs. It must be 
the triangular lattice configuration. So any compact configuration 
is formed by putting clusters of three small discs into some 
of the holes in the triangular lattice configuration of large discs. 
(Equivalently, we can form a compact configuration by removing 
some subset of the clusters of three small discs from figure \ref{fig101}.)

\medskip

\no {\bf Possible packings for $r=c_4$:} 

There is a unique sequence for small discs and 
four possible sequences for large discs:
\bea
&r : 1111 \nonumber\\
&1 : r1r1111, \quad r11r111, \quad r1r1r1r1, \quad or \quad 111111
\eea

The sequence for small discs implies that every small disc is surrounded 
by four large discs whose vertices form a square with sides of length $2$.
In the configuration in figure \ref{fig414} these squares tile the plane.
However, there are many compact packings in which they do not tile the 
plane. We will show there is a one to one correspondence between compact
packings with $r=c_4$ and tilings of the plane using  squares with sides of
length $2$ and equilateral triangles with sides of length $2$. 
An example of a packing and the corresponding tiling is shown in figure
\ref{fig414a}.

Given a packing, we draw line segments between the centers of any two 
large discs that are tangent. These line segments, which all have length
$2$, divide the plane into squares and equilateral triangles. To see this,
consider a large disc and the possible sequences around it. For the 
sequence $r1r1111$, the large disc is a vertex of two adjacent squares 
and three triangles. For the sequence $r11r111$, the large disc is again
a vertex of two squares and three triangles, but now the two squares
do not share a side. For the sequence $r1r1r1r1$, the large disc is 
a vertex of four squares. And for the sequence $111111$ the large disc is 
a vertex of six triangles. 
So in each case the center of the large disc is a vertex of squares and 
triangles that fit together to tile the space around the large disc. 

Thus the packing can be used to construct a tiling
by squares and triangles. Conversely, given such a tiling, we can construct
a compact packing by putting a large disc at the vertices of the squares and 
triangles and then putting a small disc at the center of each square. 

\medskip

\no {\bf Possible packings for $r=c_7$:} 

The sequence for small discs is unique and there are six sequences for 
larges discs:
\bea
 & r : 111r \nonumber\\
 & 1 : rr1r1111(8), \quad rr11r111(8), \quad rrr1r1r1r1(10), \quad \nonumber\\
 & rr1rr1r1r1(10), \quad rr1r1rr1r1(10) \quad or \quad 111111(6)
\eea
The third of the possibilities for a large disc implies that there is 
a small disc that is tangent to two other small discs. This contradicts 
the sequence for small discs. So the third possibility for large discs does
not occur. 

The sequence for small discs implies that each small disc is tangent to 
one other small disc and this pair is surrounded by four large discs. 
The centers of these four large discs are the vertices of a rhombus
whose sides are of length 2 and whose acute angle is
$2 \cos^{-1}({-1 + \sqrt{17} \over 4})$. Small discs can only appear 
as pairs inside such rhombi. In the packing shown in figure \ref{fig280} 
these rhombi tile the plane. 

There are many possible compact packings for $r=c_7$. We claim there is 
a one to one correspondence between compact packings and tilings of the 
plane using these rhombi and equilateral triangles with sides of length $2$. 
An example of a packing and the corresponding tiling is shown in figure
\ref{fig280a}.
Given a packing we construct the tiling by drawing line segments between 
the centers of any pair of large discs that are tangent to each other.
To see that this produces a tiling of the plane by rhombi and triangles,
we can argue as we did in the case of $r=c_4$. We consider the 
possible sequences for a large disc and see that in each case the 
center of the disc is a vertex of rhombi and triangles that fit together
to tile the space around the disc.
Given a tiling we construct a packing by putting a large disc at each 
vertex of each triangle and each rhombi. Then we put a pair of small 
discs inside each rhombi. 

\medskip

\no {\bf Possible packings for $r=c_2$:} 

The sequence for small discs is unique, and there are three possible sequences
for large discs:
\bea
& r : 111rr \nonumber\\
& 1 : rr1r111, \quad r11r11r \quad or \quad 111111 
\eea

In addition to the packing in figure \ref{fig545}, we give a second 
example in figure \ref{fig545t}.
As in the previous two cases we will show that the possible compact
packings are equivalent to a tiling problem. 
In both of the figures the small discs all appear in clusters of three. 
We will argue that this property is true in any compact packing.
The sequence for small discs implies that every small disc is tangent to 
two other small discs which are tangent to each other. 
The sequence for small discs applied to these three small discs 
then implies that they are surrounded by six large discs. The vertices
of these six large disc form a six-sided polygon with sides of length
$2$. We have drawn one in figure \ref{fig545}. 

We construct the tilings in the same way as in the previous two 
cases. We draw a line segment between the centers of any two discs of 
radius $2$ that are tangent. These lines form equilateral triangles 
in addition to the six-sided polygons. 
The same arguments used in the previous cases show that given a packing, 
the six-sided polygons and triangles tile the plane and conversely, 
given such a tiling we can construct a compact packing. 

\medskip

\no {\bf Possible packings for $r=c_1$:} 

The possible sequences are 
\bea
& r : 1111r \nonumber \\ 
& 1 : r1r1r1r \quad or \quad 111111
\eea
Start with a small disc. The unique sequence for small discs 
implies it touches one other small disc. 
By the sequence for small discs, this pair of tangent small discs is 
surrounded by six large discs. The centers of these six large discs
are the vertices of a flattened hexagon. 
(An example of the flattened hexagon is shown in figure \ref{fig637}.)
We will show that these flattened hexagons tile the plane.

The flattened hexagon has two different angles in it
which we call $\theta_1$ and $\theta_2$ with $\theta_1<\theta_2$.
We have 4 $\theta_1$'s and 2 $\theta_2$'s in the flattened hexagon. 
The sequence for large discs implies that at every large disc we must have 
two $\theta_1$'s and one $\theta_2$. It follows that the 
flattened hexagons tile the plane. 

To determine the possible tilings, we consider two cases. 
The first case is that all the flattened hexagons
have the same orientation. The second case is that they do not. 
It is easy to see that in the first case the configuration must be 
the one shown in figure \ref{fig637}.
In the second case there must be two flattened hexagons
that share a side but have different orientations. So we have two
flattened hexagons which fit together as shown by the two bold hexagons in 
figure \ref{fig637a}. Using the constraint that each large disc has
two $\theta_1$'s and one $\theta_2$, it follows that all the shaded 
hexagons must be present, extending to infinity in both directions. 

Next we consider what can happen above and below this shaded layer. 
Consider adding a flattened hexagon with two sides adjacent to the 
shaded layer. There are two possible orientations, but as soon as 
we add one flattened hexagon, the orientation of all the other ones
touching this boundary of the shaded region must be the same. 
In the figure we have added one above the shaded region with the 
same orientation as the shaded flattened hexagons along this upper boundary.
And below the shaded region we have added one whose orientation is different
from the shaded flattened hexagons along the lower boundary.
It follows that the configuration must consist of such layers with 
the freedom to chose an orientation for each layer. 

\medskip 

\no {\bf Possible packings for $r=c_3$:} 

The sequence for small discs is unique, and there are five possible 
sequences for large discs:
\bea
& r : 1r1r1 \nonumber\\
& 1 : rrr1111, \quad rrrrr1r1, \quad rrrr1rr1, \quad rrr1rrr1, 
\quad or \quad 111111
\eea
We will show that there are two types of compact configurations for 
this value of $r$. Figure \ref{fig533a} gives an example from the first 
class. In general, a configuration in this class is formed of layers of 
large discs and layers of small discs (which oscillate up and down 
slightly as in the figure). The only constraint on the sequence of 
layers is that we cannot have two adjacent layers of small discs. 
An example from the second class is shown in figure \ref{fig533b}. 
These configurations are formed by layers made up of large discs alternating
with two tangent small discs aligned vertically. (For example, the 
shaded discs in figure \ref{fig533b} consist of three layers.)  
Each layer has two possible locations relative to the layer above (or below) 
it. Note that the configuration in figure \ref{fig533} belongs to 
both classes and is the only configuration that does.

We start by showing that the sequence $1:rrrrr1r1$ does not occur. 
Consider the small disc represented
by the last $r$. It touches the large disc at the center of the sequence, 
and there are two other larges discs tangent to the small disc and
the large disc at the center. So the sequence around the small disc
must contain $111$ which is not possible. So the sequence 
$1:rrrrr1r1$ does not occur. 

We now divide into two cases. The first case is that at least one
large disc has the sequence $1:rrr1111$. The second case is that no 
large disc has this sequence. We will show that in the first case the 
configuration is of the first type described above, and in the second case
it is of the second type. 

In figure \ref{fig533a} the large disc that is assumed to have the sequence
$1:rrr1111$ is labeled $a$. It is surrounded by $b,c,d,e,f,g,h$. 
Using the possible sequences, we now have 
\bea 
e & \implies i,j,k,l \nonumber\\ 
h & \implies m,n,o,p \nonumber\\ 
c & \implies q,r 
\eea 
The large discs $e$ and $h$ are now surrounded by the same sequence that
$a$ is. Repeating this argument, we conclude that all of the shaded discs
are present. We now consider what can happen above and below the shaded 
discs. Consider two large shaded discs that are tangent to each other and 
at the edge of the shaded discs. They must both be tangent to a disc that
is not shaded. It is either small or large. 
The case of a small disc is shown by $s$ above the shaded discs, 
and the case of a large disc is shown by $t$ below the shaded discs.
It is now straightforward to show that $s$ forces the layer of small unshaded
discs just above the shaded discs as well as the layer of large unshaded
discs at the top edge of the figure. And the small disc $t$ is easily shown
to imply the layer of large unshaded discs just below the shaded discs. 

For the second case we refer to figure \ref{fig533b}. In the second 
case we are assuming that the sequence $1:rrr1111$ does not occur. 
It is easy to see that this implies that if the configuration 
contains at least one small disc, then the sequence $1:111111$ cannot
occur. So the only possible sequences around a large disc are 
$rrrr1rr1$ and $rrr1rrr1$. If the sequence $rrrr1rr1$ never appears, 
then it is easy to show that the configuration must be that shown in 
figure \ref{fig533}. So we now assume there is at least one large disc
with the sequence $rrrr1rr1$. In figure \ref{fig533b} this disc is labeled
$a$, and the discs surrounding it $b,c,d,e,f,g,h,i$.  We now have 
\bea 
b & \implies j,k \nonumber\\ 
e & \implies l,m \nonumber\\ 
g & \implies n,o \nonumber\\ 
h & \implies p \nonumber\\ 
o & \implies q,r,s,t \nonumber\\ 
c & \implies u \nonumber\\ 
u & \implies v,w,x,y \nonumber
\eea 
This shows that the large discs $o$ and $u$ have the same sequence around
them as the large disc $a$. So we can repeat this argument and conclude that 
all of the shaded discs must be present. 

We now ask what can happen above and below the shaded discs in figure 
\ref{fig533b}. The large disc $i$ must touch one more large disc and one more
small disc in addition to the shaded discs. There are two possible 
locations for the additional small disc. 

One possibility is given by the disc $z$ just above the shaded region, and 
the other possibility by the disc $z$ just below the shaded region.
It is now easy to show that the presence 
of the $z$ disc just above the shaded discs forces all of the layer of 
unshaded discs above the shaded ones. 
And the presence of $z$ as shown below the shaded discs 
forces the layer of unshaded discs below the shaded ones.
We conclude that the configuration in the second case is made up layers. 
(In the figure there are three layers with shaded discs and single layers
above and below them that are unshaded.) Each layer has two possible locations
relative to the layer below (or above) it. 

\section*{Acknowledgements}
This work was supported by the National Science Foundation (DMS-0201566).

\end{document}